\documentclass[12pt,a4paper,reqno]{article}
\def\hybrid{\topmargin 0pt      \oddsidemargin 0pt
	\headheight 0pt \headsep 0pt
	\textwidth 160true mm       
	\textheight 231true mm         
	\marginparwidth 0.0in
	\parskip 0pt plus 1pt   \jot = 1.5ex}
\hybrid

\usepackage{amssymb}
\usepackage{amsmath}
\usepackage{amsthm}

\begin{document}
\title{Invariant quantization in one and two parameters on semisimple
coadjoint orbits of simple Lie groups}

\author
{J.~Donin\\
{\normalsize Department of Mathematics, Bar-Ilan University,}\\
{\normalsize 52900 Ramat-Gan, Israel,}\\
D.~Gurevich,\\
{\normalsize ISTV, Universit\'e de Valenciennes,}\\
{\normalsize 59304 Valenciennes, France,}\\
S.~Shnider\\
{\normalsize Department of Mathematics, Bar-Ilan University,}\\
{\normalsize 52900 Ramat-Gan, Israel}
}
\date{}

\renewcommand{\(}{\left(}
\renewcommand{\)}{\right)}
\newcommand{\g}{\mathfrak{g}}
\newcommand{\hh}{\mathfrak{h}}
\newcommand{\mm}{\mathfrak{m}}
\renewcommand{\aa}{\alpha}
\newcommand{\ff}{\varphi}
\newcommand{\ssll}{\mathfrak{sl}}
\newcommand{\Lie}{\mathfrak}
\newcommand{\K}{\mathbb{K}}
\newcommand{\tb}{\mathbf{t}}
\newcommand{\Uh}{{U}_h (\g)}
\newcommand{\U}{\overline{U}_q (\g)}
\newcommand{\Utwo}{\overline{U}_q(\ssll (2))}
\newcommand{\Uthree}{\overline{U}_q(\ssll (3))}
\newcommand{\UU}{\mathcal{U}}
\newcommand{\A}{\mathcal{A}}
\newcommand{\B}{\mathcal{B}}
\newcommand{\C}{\mathbb{C}}
\newcommand{\D}{\mathcal{D}}
\newcommand{\R}{\mathbb{R}}
\newcommand{\F}{\mathcal{F}}
\renewcommand{\L}{\mathcal{L}}
\newcommand{\M}{\mathcal{M}}
\newcommand{\Z}{\mathcal{Z}}
\newcommand{\X}{\mathcal{X}}
\newcommand{\Y}{\mathcal{Y}}
\newcommand{\I}{\mathcal{I}}
\newcommand{\ot}{\otimes}
\newcommand{\<}{\langle}
\renewcommand{\>}{\rangle}
\renewcommand{\[}{[\![}
\renewcommand{\]}{]\!]}
\newcommand{\half}{\frac{1}{2}}
\newcommand{\hlf}[1]{\frac{#1}{2}}
\newcommand{\third}{\frac{1}{3}}
\newcommand{\thrd}[1]{\frac{#1}{3}}
\newcommand{\qqquad}{\quad\quad\quad}
\newcommand{\qqqquad}{\quad\quad\quad\quad}
\newcommand{\imp}{\Rightarrow}
\newcommand{\impl}{\Longrightarrow}
\newcommand{\ad}{\operatorname{ad}}
\newcommand{\id}{\operatorname{id}}
\newcommand{\Tr}{\operatorname{Tr}}
\newcommand{\End}{\operatorname{End}}
\newcommand{\bL}{\overline{L}}
\newcommand{\bP}{\overline{P}}
\newcommand{\vep}{\varepsilon}
\renewcommand{\Im}{\operatorname{Im}}
\newcommand{\Ker}{\operatorname{Ker}}
\newcommand{\sht}{\operatorname{ht}}
\renewcommand{\ll}{\lambda}

\newcommand{\bsig}{\bar{\sigma}_h}
\newcommand{\Ga}{\Gamma}
\newcommand{\bO}{\overline\Omega_\Gamma}

\newcommand{\QED}{\hspace{0.2in}\vrule width 6pt height 6pt depth 0pt 
\vspace{0.1in}}
\hfuzz1pc 
\vfuzz1.2pt 

\theoremstyle{plain} 
\newtheorem{thm}{Theorem}[section]
\newtheorem{cor}{Corollary}
\newtheorem{lemma}{Lemma}[section]
\newtheorem{propn}{Proposition}[section]
\newtheorem{axiom}{Axiom}

\theoremstyle{definition}
\newtheorem{defn}{Definition}[section]
\newtheorem{example}[defn]{Example}
\newtheorem{examples}[defn]{Examples}

\theoremstyle{remark}
\newtheorem{rem}{Remark}[section]
\newtheorem{conjecture}{Conjecture}
\renewcommand{\theconjecture}{}

\newcommand{\thmref}[1]{Theorem~\ref{#1}}
\newcommand{\secref}[1]{Section~\ref{#1}}
\newcommand{\lemref}[1]{Lemma~\ref{#1}}
\newcommand{\propref}[1]{Proposition~\ref{#1}}

\newcommand{\be}[1]{\begin{eqnarray#1}}
\newcommand{\ee}[1]{\end{eqnarray#1}}

\newcommand{\De}{\Delta}
\newcommand{\cA}{\mathcal{C}}
\renewcommand{\t}{\otimes}
\newcommand{\wt}{\widetilde}

\newcommand{\pb}{\beta^\prime}
\newcommand{\bpb}{\overline{\beta^\prime}}
\newcommand{\bbe}{{\bar\beta}}
\newcommand{\bga}{{\bar\gamma}}
\newcommand{\Ea}{E_{\aa}}
\newcommand{\Ema}{E_{-\aa}}
\newcommand{\rc}{\Omega^+\setminus\Omega_\Ga}
\newcommand{\baa}{{\bar\aa}}

\numberwithin{equation}{section}
\maketitle

\begin{abstract}
Let $\A$ be the function algebra on a semisimple orbit, $M$,
in the coadjoint representation of a simple Lie group, $G$, with the Lie
algebra $\g$. We study one and two parameter quantizations of $\A$,
$\A_h$, $\A_{t,h}$, such that the multiplication on the quantized algebra 
is invariant under action of the Drinfeld-Jimbo
quantum group, $U_h(\g)$. In particular, the algebra $\A_{t,h}$
specializes at $h=0$ to a $U(\g)$, or $G$, invariant quantization, $A_{t,0}$.

We prove that the Poisson bracket corresponding to $\A_h$ must be
the sum of the so-called $r$-matrix and an invariant bracket. We  
classify such brackets for all semisimple orbits, $M$, and 
show that they form a $\dim H^2(M)$ parameter family, 
then we construct their quantizations. 

A two parameter (or double) quantization, $\A_{t,h}$, corresponds 
to a pair of compatible Poisson brackets: the first is as described above 
and the second is the Kirillov-Kostant-Souriau bracket on $M$. 
Not all semisimple orbits admit a compatible pair of Poisson brackets.
 We classify the semisimple orbits for which such pairs
 exist and construct the corresponding two parameter
quantization of these pairs in some of the cases.
\end{abstract}

\section{Introduction}

Passing from  classical mechanics to quantum mechanics involves
replacing the commutative function algebra, $\A$, of classical observables
on the appropriate phase space, $M$, with a noncommutative (deformed) algebra,
$\A_t$, of quantum observables (see \cite{BFFLS} where the 
deformation quantization scheme is developed). 
The algebra $\A$ is a Poisson algebra and the product in $\A_t$ 
is given by a power series in the formal parameter $t$ with leading term 
the original commutative product and with leading term in the
 commutator   given by the Poisson  bracket.
If the classical system is invariant under a 
Lie group of symmetries $G$, the associated quantum system often retains 
the group of symmetries. In particular, the algebra $\A_t$ is often 
invariant under the action of $G$, or under the action of its universal 
enveloping algebra $U(\g)$.

Modern field-theoretical models and, in particular, the problem of 
incorporating gravity into a quantum field theory led to the requirement
of deforming  (quantizing) the group symmetry and the phase space themselves.
This is one of the reasons for the interest in quantum groups. 
The quantum group, $\Uh$, defined by Drinfeld and Jimbo is a deformation of 
$U(\g)$ as a Hopf algebra.
The quantization of the phase space and its  symmetry group corresponds to
a $\Uh$ invariant deformation of the  algebra $\A_t$, which leads us
to the problem of two parameter (or double) quantization, $\A_{t,h}$, 
of the function algebra for a $U(\g)$ invariant Poisson structure.
In other words, the problem of two parameter quantization appears if
we want to quantize the Poisson bracket in such a way that 
multiplication in the quantized algebra is invariant under the 
quantum group action.

In the present paper we investigate the problem of one and two parameter
invariant quantizations of the Poisson function algebra on
a semisimple orbit in the coajoint representation of a simple Lie group. 

In general, if $M$ is a manifold on which a semisimple Lie group
$G$ acts, it is not true that there exists even a one parameter $\Uh$
invariant quantization of the function algebra, $\A_h$.
In \cite{DGM} it is proven that such a quantization exists if
$M=G/H$ and the Lie algebra of $H$ contains a maximal nilpotent
subalgebra. A $\Uh$ invariant quantization 
of the algebra of holomorphic sections of a line bundle over a flag manifold
is constructed by similar methods in \cite{DG1}.
In all these cases the Poisson bracket which is  quantized
is the $r$-matrix bracket, $r_M$, on $M$. This bracket is determined
by the bivector field $(\rho\t\rho)(r)$ for $r\in\wedge^2\g$,
the Drinfeld-Jimbo  $r$-matrix (of the form (\ref{rm})), and 
$\rho:\g\to{\rm Vect}(M)$ is the mapping defined by the action
of $G$ on $M$.

The Sklyanin-Drinfeld (SD)  bracket on $G$ is determined 
by the difference of two bivector fields on $G$
which are the right and left invariant extensions of the Drinfeld-Jimbo
classical $r$-matrix. For details see \cite{LW}. In \cite{DG2} it is shown 
that a one parameter quantization of the SD bracket exists
for all semisimple orbits in coadjoint representation of $G$.

In the present paper we show that for
the semisimple orbits, $M$, the SD bracket is only one from $\dim H^2(M)$
parameter family of Poisson brackets admitting $\Uh$ invariant
quantization.

For symmetric spaces $M$, the $r$-matrix bracket
satisfies the Jacobi identity and hence defines a Poisson bracket.
The existence of a one parameter $\Uh$ invariant quantization with the $r$-matrix Poisson 
bracket is proven in \cite{DS1}. That paper also proves that
when $M$ is a hermitian symmetric space there exists a two parameter
$\Uh$ invariant quantization.
 The hermitian symmetric spaces form a subclass of semisimple orbits
in $\g^*$. These orbits have a very interesting property: the one-sided
invariant components of the
Sklyanin-Drinfeld bracket being reduced on such an orbit become
Poisson brackets separately. 
More precisely, one of these components being reduced on $M$ is
the $r$-matrix
bracket and the other one becomes the Kirillov-Kostant-Souriau (KKS) 
bracket which is roughly speaking the restriction to the orbit of the 
Lie bracket in $\g$. 
(see \cite{KRR}, \cite{DG1}). These brackets are obviously
compatible (i.e. their Schouten bracket vanishes). So,  
we get a Poisson pencil which is the set of linear combinations
of KKS and $r$-matrix Poisson brackets. 
The two parameter quantization, $\A_{t,h}$, for hermitian symmetric
spaces constructed in \cite{DS1} is just a quantization of such a pencil.
In particular,
$\A_{t,0}$ is a $U(\g)$ invariant quantization of the KKS bracket on $M$.

The orbits in $\g^*$ on which the $r$-matrix bracket is Poisson 
have been classified in \cite{GP}.  In particular, the only such
semisimple orbits  are symmetric spaces.

It turns out that there exist two parameter $\Uh$ invariant quantizations
on some $G$-manifolds where the $r$-matrix bracket is not Poisson (does not
satisfy the Jacobi identity) but one can add a $G$-invariant bracket 
to the $r$-matrix bracket and get a Poisson bracket.
For example, in \cite{Do} it is shown that for $\g=sl(n)$
there is a two parameter $\Uh$ invariant
family $(S\g)_{t,h}$, where $S\g$ is the symmetric algebra
of $\g$ which can be considered as the function (polynomial)
algebra on $\g^*$. It is proven that this family can be restricted
to give a two parameter quantization on any semisimple orbit of
maximal dimension (on which the $r$-matrix bracket is not Poisson).

In the present paper we prove the existence of  a two parameter $\Uh$ invariant 
quantization of the function algebra  for some 
non-symmetric semisimple orbits in $\g^*$. Similar to
the case  of symmetric orbits, the quasiclassical (infinitesimal)
 term of this quantization is 
a Poisson pencil generated by the KKS bracket and another Poisson bracket 
which must be the sum of $r$-matrix and an $U(\g)$ invariant brackets.
In particular, we shall see that such pencils exist for all semisimple 
orbits in case $\g=sl(n)$.  Note, that for non-symmetric orbits both 
the $r$-matrix bracket and the SD Poisson bracket are not compatible 
with the KKS bracket. We  give a complete 
classification  of the orbits admitting such a Poisson pencil.
Moreover, we classify all such pencils and construct 
the deformation quantization of some of them.

We now describe the content of the paper in more detail.
Let $G$ be a simple, connected, complex Lie group, $\g$ its Lie algebra,
The Drinfeld-Jimbo quantum group $\Uh$ can be considered
as algebra $U(\g)[[h]]$ with undeformed multiplication but
 deformed noncocommutative comultiplication $\De_h$.

Let $M$ be a $G$-homogeneous complex manifold.
It is easy to show that $M$ is isomorphic to a semisimple
orbit of $G$ in the coadjoint representation $\g^*$ if and
only if the stabilizer, $G_o$, of a point ${\bf o}\in M$ is a Levi subgroup 
in $G$ (see Section 3 for definition). The symplectic structure on $M$
is not  unique, but each symplectic structure
on $M$ arises from an isomorphism of $M$ with an orbit, $O_\lambda$, 
for some  semisimple element $\lambda\in\g^*$. 
On $O_\lambda$ there is
the Kirillov-Kostant-Souriau (KKS) Poisson bracket, $v_\lambda$,
whose action on the restriction of linear functions to the orbit 
is given by the Lie bracket in $\g$.
Each symplectic structure on $M$ is induced from the KKS bracket
by an isomorphism of $M$ onto a semisimple orbit
(see Section 3).

The first problem we consider is that of quantizing  
the algebra $\A$ of polynomial (or holomorphic) functions on $M$, such 
that the quantized algebra $\A_h$ has a $\Uh$ invariant multiplication
$$\mu_h=\sum_{i=0}^\infty h^i\mu_i,$$ 
where $\mu_0$ is the initial multiplication of functions and the 
$\mu_i$ for $i\geq 1$ are bidifferential operators.
Invariance means that $\mu_h$ satisfies the property
\be{*}
\mu_h\De_h(x)(a\t b)=x\mu_h(a\t b)\qqquad {\rm for} \quad
x\in\Uh,\quad a,b\in \A.
\ee{*} 

The second problem we consider is the existence of a two parameter 
$\Uh$ invariant quantization, $\A_{t,h}$, of $\A$ such that the 
one parameter family $\A_{t,0}$ is a $U(\g)$ invariant quantization 
of the KKS bracket on $M$.

The infinitesimal of the one parameter quantization $\A_h$ is
a Poisson bracket, whereas the infinitesimal of the two parameter quantization
is a pencil of two compatible Poisson brackets. 

In Section 2 we recall some facts on the Drinfeld monoidal categories
of quantum group representations, using in our construction of quantization.
In addition, we show that the Poisson bracket, $p(a,b)=\mu_1(a,b)-\mu_1(b,a)$,
corresponding to the $\Uh$ invariant quantization $\A_h$ 
must be of a special form. Namely, let $r\in\wedge^2\g$ be 
the Drinfeld-Jimbo classical $r$-matrix. The Schouten
bracket $\[r,r\]\in\wedge^3\g$ is the invariant element $\ff$, which is
unique up to a factor. Denote by $r_M$ the bracket on $M$
determined by the bivector field $(\rho\t\rho)(r)$ where 
$\rho:\g\to{\rm Vect}(M)$ is the mapping defined by the action
of $G$ on $M$. We call $r_M$ an $r$-matrix bracket.
Put also $\ff_M=\rho^{\t 3}\ff$.  Then, $p(a,b)$ has the form
\be{}\label{e1}
p(a,b)=r_M(a,b)+f(a,b), \qqquad a,b\in \A,
\ee{}
where $f(a,b)$ is a $U(\g)$ invariant bracket
with the Schouten bracket 
\be{}\label{e2}
\[f,f\]=-\ff_M. 
\ee{}
Note that $r$-matrix
bracket is compatible with any invariant bracket,  
i.e. $\[f,r_M\]=0$. 

Similarly, the two parameter quantization $\A_{t,h}$ corresponds to
a pair of compatible Poisson brackets, $(p,v_\lambda),$
where $p=r_M+f$ of the form (\ref{e1}) and $v_\lambda$ is
the KKS bracket. Since $r_M$ is compatible with
the invariant bracket $v_\lambda$, compatibility of $p$ with $v_\lambda$
is equivalent to condition 
\be{}\label{e3}
\[f,v_\lambda\]=0. 
\ee{}

In Section 3 we give a classification of all invariant brackets
$f$ satisfying condition (\ref{e2}) for all $M$ isomorphic to
semisimple orbits. 
We show that such brackets form a $\dim H^2(M)$ parameter family.

In the same section we show that a semisimple orbit may not have
any invariant brackets satisfying the
conditions (\ref{e2}) and (\ref{e3}).
 We call an orbit a ``good orbit'', if such a bracket exists and then
give a classification of all good orbits. Namely,
if $\g$ is of type $A_n$, all semisimple orbits are good.
All orbits which are symmetric spaces are good. In cases
$D_n$, $E_6$ there are good orbits which are not symmetric spaces.
Moreover, for the good orbits the brackets satisfying the conditions
(\ref{e2}) and (\ref{e3}) form a one parameter family.
In fact, the property of an orbit being  good  depends only on its structure
as a homogeneous manifold, not on the symplectic structure. The 
dependence is on the Lie subalgera of the stabilizer subgroup. So, 
if an orbit $M$ is good, 
then any orbit isomorphic to $M$ as a homogeneous manifold will be good.

In Section 4 we consider cohomologies of the complex of invariant polyvector
fields on $M$ with differential given by the Schouten bracket with
the bivector $f$ satisfying (\ref{e2}). We show that for almost all
$f$ these cohomologies coincide with the usual de Rham cohomologies
of the manifold $M$ and then use this fact in Section 5 to prove the existence
of an invariant quantization. In the proof we use methods of 
\cite{DS1} and \cite{DS2}.
Using the same methods, we also construct the two parameter quantization
for good orbits in cases $D_n$, $E_6$ and brackets satisfying
(\ref{e2}) and (\ref{e3}). 
For the case $A_n$ some additional arguments are required.
(See \cite{Do} where using another method the existence of
two parameter quantization for maximal orbits is proven.)

In conclusion we make two remarks.

{\it Remark 1}. In the paper $G$ is supposed to be a complex Lie group. 
However, one can consider the situation when $G$ is a real
simple Lie group with Lie algebra $\g_\R$ and $M$
is a semisimple orbit of $G$ in $\g^*_\R$. 
In this case we take $\A=C^\infty(M)$,
the complex-valued smooth functions. Let $\g$ be the complexification of $\g_\R$.
It is clear that Lie algebra $\g$ and algebra $U(\g)$ act on $C^\infty(M)$.
Since all our results are formulated in terms of $\g$ action, they are
valid in the real case as well (see also \cite{DS1}).

{\it Remark 2}. 
The deformation quantization can be considered as the first
 step of a quantization procedure whose the second step is 
a representation of the quantized 
algebra, $\A_{t,h}$, as an operator algebra in a linear space.
For some symmetric orbits in $sl(n)^*$ such a representation has been 
given in \cite{DGR} and
\cite{DGK} but the method of \cite{DGK} can be apparently extended to all
symmetric orbits in $\g^*$ for all simple Lie algebras $\g$.
These operator algebras have a deformed $U_h(\g)$ invariant
trace which, however, is not symmetric.

\section{Poisson brackets associated with $\Uh$ invariant\\ quantization}
We recall some facts about the Drinfeld algebras and the monoidal
categories determined by them. They will be used, in particular, in our
construction of the quantization.

Let $A$ be a commutative algebra with unit, $B$ a unitary
$A$-algebra. The category of representations of $B$ in $A$-modules,
i.e. the category of $B$-modules, will be a monoidal category 
if the algebra $B$ is equipped with an algebra morphism, $\De: B\to B\ot_A B$, 
 called  comultiplication,  and an invertible element
$\Phi\in B^{\t 3}$ such that $\De$ and $\Phi$ 
satisfy the conditions (see \cite{Dr1})
\be{}
&(id\t\De)(\De(b))\cdot\Phi=\Phi\cdot (\De\t id)(\De(b)),\ \ b\in B, 
\label{d1}\\
&(id^{\t 2}\t \De)(\Phi)\cdot(\De\t id^{\t 2})(\Phi)=
(1\t\Phi)\cdot(id\t\De\t id)(\Phi)\cdot(\Phi\t 1).  \label{d2}       
\ee{}
Define a tensor product functor
for $\cA$ the category of
$B$ modules,  denoted  $\ot_{\cA}$  or simply $\t$ when 
there can be no confusion, in
the following way: given $B$-modules $M,N$,
$M\t_{\cA} N=M\t_A N$ as an $A$-module. The
action of $B$ is defined by
 $$b(m\t n)=(\De b)(m\t n)= b_1m\t b_2n\quad\mbox{ where }\De b=b_1\t b_2,$$
using the Sweedler convention of an implicit summation over an index.
The element $\Phi=\Phi_1\t\Phi_2\t\Phi_3$ defines the
 associativity constraint,
$$a_{M,N,P}:(M\t N)\t P\to M\t(N\t P),\,\, 
a_{M,N,P}((m\t n)\t p)= \Phi_1m\t(\Phi_2n\t\Phi_3p).$$
Again the summation in the expression for $\Phi$ is understood.
By virtue of (\ref{d1}) $\Phi$ induces an isomorphism of $B$-modules, and by
virtue of
(\ref{d2}) the pentagon identity for monoidal categories holds. We call the
triple
$(B,\De,\Phi)$ a Drinfeld algebra.
The definition is somewhat non-standard in that
 we do not require the existence of an antipode.
The category $\cA$ of $B$-modules for $B$ a 
Drinfeld algebra  becomes a monoidal category. When
it becomes necessary to be more explicit  we shall denote 
${\cA}(B,\De,\Phi)$. 

Let $(B,\De,\Phi)$ be a Drinfeld algebra and $F\in B^{\t 2}$
an invertible element. Put 
\be{}\label{f1}
&\wt{\De}(b)=F\De(b)F^{-1},\ \ b\in B,  \\
&\wt{\Phi}=(1\t F)\cdot(id\t\De)(F)\cdot\Phi\cdot(\De\t id)(F^{-1})
\cdot(F\t 1)^{-1}.           \label{f2}
\ee{}
Then $\wt{\De}$ and $\wt{\Phi}$ satisfy (\ref{d1}) and (\ref{d2}),
therefore the triple $(B,\wt{\De},\wt{\Phi})$ also becomes a Drinfeld
algebra which has an equivalent monoidal category of modules, 
$\wt{\cA}(B,\wt{\De},\wt{\Phi})$. Note that the equivalent categories $\cA$
and $\wt{\cA}$ consist of the same objects as $B$-modules,
and the tensor products of two objects are isomorphic as $A$-modules.
 The equivalence $\cA\to \wt{\cA}$ is given by the pair $(Id,F)$,
where $Id:\cA\to \wt{\cA}$ is the identity functor of the categories
(considered without the monoidal structures, but only as categories
of $B$-modules), and $F:M\t_{\cA} N\to M\t_{\wt{\cA}} N$ is defined by 
$m\t n\mapsto F_1m\t F_2n$ where $F_1\t F_2=F$. 

Assume $\A$ is a $B$-module with a multiplication $\mu:\A\t_A \A\to \A$
which is a homomorphism of $A$-modules. We say that $\mu$
is $\De$ invariant if 
\be{}\label{finv}
b\mu(x\otimes y)=\mu\De(b)(x\t y)\ \ \ \mbox{for}\ b\in B,\ x,y\in \A,
\ee{}
and $\Phi$ associative, if
\be{}\label{fass}
\mu(\Phi_1x\otimes \mu(\Phi_2 y\t \Phi_3 z)))=\mu(\mu(x\t y)\t z)
\ \ \ \mbox{ for } x,y,z\in \A. 
\ee{}

Note, that a $B$-module $\A$ equipped with $\De$ invariant and
$\Phi$ associative multiplication is an associative algebra in
the monoidal category ${\cA}(B,\De,\Phi)$.
The multiplication $\wt{\mu}=\mu F^{-1}:M\t_A M\to M$ will be
$\wt{\Phi}$-associative and invariant in the category $\wt{\cA}$. 

We are interested in the case when $A=\C[[h]]$, $B=U(\g)[[h]]$
where $\g$ is a complex simple Lie algebra. In this case, all
tensor products over $\C[[h]]$ are completed in
$h$-adic topology.

Denote by $\ff\in\wedge^{\t 3}\g$ an invariant element (unique up to
scaling for $\g$ simple) and
by $r\in\wedge^{\t 2}\g$ the so-called Drinfeld-Jimbo $r$-matrix
of the form (\ref{rm}) such that the Schouten bracket of $r$ with itself
is equal to $\ff$:
\be{}\label{eq2.7}
\[r,r\]=\ff.
\ee{}
In \cite{Dr1}, Drinfeld proved the following (see also \cite{DS2} 
for the property c)).
\begin{propn}
\label{p1.1}
1. There is an invariant element $\Phi_h\in U(\g)[[h]]^{\t 3}$ of the form 
$\Phi_h=1\t 1 \t 1+h^2\ff+\cdots$  satisfying the following
properties:

\ a) it depends on $h^2$, i.e. $\Phi_h=\Phi_{-h}$; 

\ b) it satisfies the equations (\ref{d1}) (i.e. invariant) and (\ref{d2}) 
with the usual $\De$ arising from $U(\g)$;

\ c) it is invariant under the Cartan involution $\theta$;

\ d) $\Phi^{-1}_h=\Phi^{321}_h$, where 
$\Phi^{321}=\Phi_3\t\Phi_2\t\Phi_1$ for $\Phi=\Phi_1\t\Phi_2\t\Phi_3$;

2. There is an element $F_h\in U(\g)[[h]]^{\t 2}$ of the form 
$F_h=1\t 1+(h/2)r+\cdots$ satisfying the 
equation (\ref{f2}) with the usual $\De$ and with 
$\wt{\Phi}=1\t 1\t 1$.
\end{propn}

This proposition implies that there
are two nontrivial Drinfeld algebras: the first, $(U(\g)[[h]],\De,\Phi_h)$
with the usual comultiplication and $\Phi$ from Proposition
\ref{p1.1}, and the second, $(U(\g)[[h]],\wt{\De}, {\bf 1})$ 
where $\wt{\De}(x)=F_h\De(x)F_h^{-1}$ 
for $x\in U(\g)$. The pair $(Id,F_h)$ defines an
equivalence between the  corresponding monoidal categories 
${\cA}(U(\g)[[h]],\De,\Phi_h)$
and ${\cA}(U(\g)[[h]],\wt{\De},{\bf 1})$ 

It is clear that reduction modulo $h$
defines a functor from  either of these categories
to the category of representations of $U(\g)$ and the equivalence
just described reduces to the identity modulo $h$. 
In fact, both categories are 
$\C[[h]]$-linear extensions of the $\C$-linear
category of representations of $\g$. Ignoring the monoidal structure
the extension is a trivial one, but the associator $\Phi$ in the first case
and the  comultiplication $\wt\De$ in the second case make the
extension non-trivial from the point of view of monoidal categories.

The bialgebra $U(\g)[[h]]$ with comultiplication $\De_h=\wt{\De}$ is denoted
by $\Uh$ and is isomorphic to the Drinfeld-Jimbo quantum group (\cite{Dr1}).

Let $\A$ be a $U(\g)$ invariant commutative algebra, i.e. an algebra with 
$U(\g)$ invariant multiplication $\mu$ in sense
of (\ref{finv}). A quantization of $\A$ is
an associative algebra,  $\A_h$, which is 
isomorphic to $\A[[h]]=\A\t\C[[h]]$ (completed tensor
product) as a $\C[[h]]$-module, with multiplication in
$\A_h$ having  the form
$\mu_h=\mu+h\mu_1+o(h)$.
The Poisson bracket corresponding to the quantization is given by
$\{a,b\}=\mu_1(a,b)-\mu_1(b,a)$, $a,b\in\A$.

In general, we call a skew-symmetric bilinear form $\A\t\A\to\A$ a bracket, 
if it satisfies the Leibniz rule in either argument when the other is fixed.
 The term Poisson bracket indicates that the Jacobi identity is also true.
A bracket of the form 
\be{}\label{rmb}
\{a,b\}_r=(r_1a)(r_2b)=\mu(r_1a,r_2b) \qqquad a,b\in\A, 
\ee{}
where $r=r_1\t r_2$ (summation implicit) is the representation of 
$r$-matrix $r$ will be called an $r$-matrix bracket.

Assume, $\A_h$ is a $\Uh$ invariant quantization, i.e. the multiplicatin
$\mu_h$ is $\De_h$ invariant. 
We shall show that in this case the Poisson bracket $\{\cdot,\cdot\}$ has
a special form.  Suppose $f$ and $g$ are two brackets
on $\A$. Then we define their Schouten bracket $\[f,g\]$ as 
\be{}
\[f,g\](a,b,c)=f(g(a,b),c)+g(f(a,b),c)+{\rm cyclic\ permutations\ of}\ a,b,c.
\ee{}
Then $\[f,g\]$ is a  skew-symmetric map $\A^{\t 3}\to \A$.
We call $f$ and $g$ compatible, if $\[f,g\]=0$.

\begin{propn}\label{p1.2}
Let $\A$ be a $U(\g)$ invariant commutative algebra and
$\A_h$ a $\Uh$ invariant quantization.
Then the corresponding Poisson bracket has the form
\be{}\label{rinv}
\{a,b\}=f(a,b)-\{a,b\}_r
\ee{}
where $f(a,b)$ is a $U(\g)$ invariant bracket.

The brackets $f$ and $\{\cdot,\cdot\}_r$ are
compatible and $\[f,f\]=-\ff_\A$
where $\ff_\A(a,b,c)=(\ff_1a)(\ff_2b)(\ff_3c)$ 
and $\ff_1\t\ff_2\t\ff_3=\ff\in\wedge^3\g$ is the invariant element.
\end{propn}
\begin{proof}
The permutation $\sigma:\A_h^{\t 2}\to \A_h^{\t 2}$, $a\t b\to b\t a$,
is an equivariant operator in the category ${\cA}(U(\g)[[h]],\De,\Phi_h)$,
because $\De$ is a cocommutative comultiplication.
The equivalence of categories
${\cA}(U(\g)[[h]],\De,\Phi_h)$ and
${\cA}(U(\g)[[h]],\wt{\De},{\bf 1})$
implies that the operator $\wt{\sigma}=F\sigma F^{-1}$ on $\A^{\t 2}$
is equivariant under
action of the quantum group $\Uh$.

Suppose the multiplication in $\A_h$ has the form
$\mu_h(a,b)=ab+h\mu_1(a,b)+o(h)$.
It is easy to calculate that
\be{*}
\frac{1}{h}\mu_h(Id-\wt{\sigma})(a\t b)=
\mu_1(a,b)-\mu_1(b,a)+(r_1a)(r_2b)+O(h)=\{a,b\}+\{a,b\}_r +O(h).
\ee{*} 
But this is a $\Uh$ equivariant operator $\A_h^{\t 2}\to \A_h$.
Taking $h=0$ we obtain that the bracket
$f(a,b)=\{a,b\}+\{a,b\}_r$ must be $U(\g)$ invariant.
So, we have
$\{a,b\}=f(a,b)-\{a,b\}_r$, as required.

It is easy to check that any bracket of the form
$\{a,b\}=(X_1a)(X_2b)=\mu(X_1a,X_2b)$, for $X_1\t X_2\in \g\wedge \g$, 
 is compatible with any invariant bracket. 
In particular, an $r$-matrix bracket is compatible with $f$.
In addition, $\{\cdot,\cdot\}$ is a Poisson bracket, so its Schouten
bracket with itself is equal to zero. Using this and 
taking into account that the Schouten bracket of $r$-matrix bracket with
itself is equal to $\ff_\A$, we obtain from (\ref{rinv}) that
$\[f,f\]=-\ff_\A$.
\end{proof}

\begin{rem}\label{rem1.1}
a) It is clear that if $r$ satisfies (\ref{eq2.7}), then $-r$ satisfies
(\ref{eq2.7}), too, and we may replace $F_h$ in Proposition
\ref{p1.1} with $\wt{F_h}$ with leading terms $1\otimes 1 -(h/2)r$.
Then, instead (\ref{rinv}) we can write
$\{a,b\}=f(a,b)+\{a,b\}_r$.

b) Assume that $\A_{t,h}$ is a two parameter quantization of $\A$,
i.e. a topologically free $\C[[t,h]]$-module 
with a multiplication of the form
$\mu_{t,h}(a,b)=ab+h\mu_1(a,b)+t\mu^\prime_1(a,b)+o(t,h)$.
Assume, that $\A_{t,h}$ is $\Uh$ invariant, so that $\A_{t,0}$ is $U(\g)$
invariant. Then there are two compatible Poisson brackets 
corresponding to such a quantization:
the bracket $\mu_1(a,b)-\mu_1(b,a)$ of the form (\ref{rinv})
and the $U(\g)$ invariant bracket 
$v(a,b)=\mu^\prime_1(a,b)-\mu^\prime_1(b,a)$. Since $v$ is invariant,
the compatibility is equivalent to $\[f,v\]=0$.

c) In view of equivalence of the categories
${\cA}(U(\g)[[h]],\De,\Phi_h)$ and
${\cA}(U(\g)[[h]],\wt{\De},{\bf 1})$,
the problem of quantizing the algebra $\A$ may be considered in the first
category. If $\A_h$ is a $\Uh$ invariant quantization with multiplication
$\mu_h$, then the multiplication 
$\bar{\mu}_h=\mu_h F_h=\mu+h\bar{\mu}_1+o(h)$
will be $U(\g)$ invariant and
$\Phi_h$ associative in sense of (\ref{fass}).
We have $\bar{\mu}_1(a,b)-\bar{\mu}_1(b,a)=f(a,b)$,
where $f$ is from (\ref{rinv}).
So, we see that the invariant bracket $f$ from (\ref{rinv}) 
with $\[f,f\]=-\ff_\A$ plays 
the role of Poisson bracket for $\Phi_h$ associative quantization.
Similarly, the two parameter quantization $\A_{t,h}$ corresponds to 
the $U(\g)$ invariant $\Phi_h$ associative quantization in
${\cA}(U(\g)[[h]],\De,\Phi_h)$ with
a pair of compatible invariant brackets $f$ and $v$ where
$\[f,f\]=-\ff_\A$ and $\[v,v\]=0$.
Working in the category ${\cA}(U(\g)[[h]],\De,\Phi_h)$
can simplify the process of quantization.
(see \cite{DS1} and Section 5).
\end{rem}

In the next section we consider the case when $\A$
is a function algebra on a semisimple orbit in the coadjoint representation
of $\g$ and we give a classification of the invariant brackets $f$ satisfying
the property $\[f,f\]=-\ff_\A$. Moreover, among such $f$ we 
distinguish those which are compatible with the KKS Poisson brackets.

\section{Pairs of brackets on semisimple orbits}
Let $\g$ be a simple complex Lie algebra, $\hh$ a fixed Cartan subalgebra.
Let $\Omega\subset \hh^*$ be  the system of roots corresponding to $\hh$.
Select a system of positive roots,  $\Omega^+$, and 
denote by $\Pi\subset\Omega$ the subset of simple roots. 
Fix an element $E_\aa\in \g$  of weight $\alpha$ for each $\alpha\in\Omega^+$ 
and choose $E_{-\aa}$ such that
$(E_\aa,E_{-\aa})=1$ for the Killing form $(\cdot,\cdot)$ on $\g$.
Then, for all pairs of roots $\aa,\beta$ such that $\aa+\beta\neq 0$
we define the numbers $N_{\aa,\beta}$ in the following way:
\be{*}
[E_\aa,E_\beta]=N_{\aa,\beta}E_{\aa+\beta} \qquad {\rm if} \quad 
\aa+\beta\in \Omega, \\
N_{\aa,\beta}=0 \qquad {\rm if} \quad \aa+\beta\notin \Omega.
\ee{*}

These numbers satisfy the following property, \cite{He}. For the roots
$\aa,\beta,\gamma$ such that $\aa+\beta+\gamma=0$ one has
\be{} \label{he}
N_{\aa,\beta}=N_{\beta,\gamma}=N_{\gamma,\aa}.
\ee{}

Let $\Ga$ be a subset of $\Pi$. Denote by $\hh^*_\Ga$
the subspace in $\hh^*$ generated by $\Ga$. 
Note, that $\hh^*=\hh^*_\Ga\oplus\hh^*_{\Pi\setminus\Ga}$,
and one can identify $\hh^*_{\Pi\setminus\Ga}$ and $\hh^*/\hh^*_\Ga$
via the  projection $\hh^*\to\hh^*/\hh^*_\Ga$.

Let $\Omega_\Ga\subset \hh^*_\Ga$ 
be the subsystem of roots
in $\Omega$ generated by $\Ga$, i.e. 
$\Omega_\Ga=\Omega\cap \hh^*_\Ga$. 
Denote by $\g_\Ga$ the subalgebra of $\g$
generated by the elements $\{E_\aa,E_{-\aa}\}$, $\aa\in\Ga$, and $\hh$.
Such a subalgebra is called the Levi subalgebra.

Let $G$ be a complex connected Lie group with Lie algebra $\g$ and  
$G_\Ga$ a subgroup with Lie algebra $\g_\Ga$. Such a subgroup
is called the Levi subgroup. It is known that
$G_\Ga$ is a connected subgroup. 
Let $M$ be a homogeneous space
of $G$ and $G_\Ga$ be the stabilizer of a point $o\in M$. We
can identify $M$ and the coset space $G/G_\Ga$. 
It is known, that such $M$ is isomorphic to a semisimple orbit in $\g^*$.
This orbit goes through an element $\lambda\in\g^*$ which is just the 
trivial extension to all of  $\g^*$ (identifying
$\g$ and $\g^*$ via the Killing form) of a map 
$\lambda:\hh_{\Pi\setminus\Ga}\to\C$ such that $\lambda(\alpha)\neq 0$
for all $\alpha\in \Pi\setminus\Ga$.
Conversely, it is easy to show that any semisimple orbit in $\g^*$ 
is isomorphic to the quotient of $G$ by a Levi subgroup.

The projection $\pi:G\to M$ induces the map
$\pi_*:\g\to T_o$ where $T_o$ is the tangent space to $M$ at the
point $o$. Since the $\ad$-action of $\g_\Ga$ on $\g$ is semisimple,
there exists an $\ad(\g_\Ga)$-invariant subspace $\mm=\mm_\Ga$ of $\g$
complementary to $\g_\Ga$, and one can identify $T_o$ and $\mm$ by means of
$\pi_*$. It is easy to see that subspace $\mm$ is uniquely defined
and has a basis formed by the elements
$E_\gamma,E_{-\gamma}$, $\gamma\in\Omega^+\setminus\Omega_\Ga$.

Let $v\in\g^{\ot m}$ be a tensor over $\g$. Using
the right and the left actions of $G$ on itself, one can associate with $v$
right and left invariant tensor fields on $G$ denoted by $v^r$ and $v^l$.

We say that a tensor field, $t$, on $G$ is right $G_\Ga$ invariant, if
$t$ is invariant under the right action of $G_\Ga$.
The $G$ equivariant diffeomorphism between $M$ and $G/G_\Ga$ implies that
any right $G_\Ga$ invariant tensor field $t$ on $G$
induces tensor field $\pi_*(t)$ on $M$. The field $\pi_*(t)$ will be
invariant on $M$ if, in addition, $t$ is left invariant on $G$, 
and any invariant
tensor field on $M$ can be obtained in such a way.
Let $v\in\g^{\ot m}$. For $v^l$ to be right $G_\Ga$
invariant it is necessary and sufficient that $v$ to be $\ad(\g_\Ga)$
invariant. Denote $\pi^r(v)=\pi_*(v^r)$ for any tensor $v$ on $\g$
and $\pi^l(v)=\pi_*(v^l)$ for any $\ad(\g_\Ga)$ invariant tensor 
$v$ on $\g$. 
Note, that tensor $\pi^r(v)$ coincides with the image of $v$ by
the map $\g^{\t m}\to{\rm Vect}(M)^{\t m}$ induced by the action map 
$\g\to{\rm Vect}(M)$.
Any $G$ invariant tensor on $M$ has the form
$\pi^l(v)$. Moreover, $v$ clearly can be uniquely choosen from $\mm^{\ot m}$. 

Denote by $\[v,w\]\in\wedge^{k+l-1}\g$ the Schouten bracket 
of the polyvectors $v\in\wedge^k\g$, $w\in\wedge^l\g$, defined
by the formula 
\be{*}
\[X_1\wedge\cdots\wedge X_k, Y_1\wedge\cdots\wedge Y_l\]=\sum
(-1)^{i+j}[X_i,Y_j]\wedge X_1\wedge\cdots \hat X_i \cdots
\hat Y_j\cdots \wedge Y_l,
\ee{*}
where $[\cdot,\cdot]$ is the bracket in $\g$.
The Schouten bracket is defined in the same way for polyvector fields on
a manifold, but instead of $[\cdot,\cdot]$ one uses the Lie bracket of
vector fields. We will use the same notation for the Schouten
bracket on manifolds.
It is easy to see that $\pi^r(\[v,w\])=\[\pi^r(v),\pi^r(w)\]$, and
the same relation is valid for $\pi^l$.

Denote by $\bO$ the image of $\Omega$ in $\hh^*_{\Pi\setminus\Ga}$
without zero.
It is clear that $\Omega_{\Pi\setminus\Ga}$ can be identified with
a subset of $\bO$ and each element
from $\bO$ is a linear combination of elements from $\Pi\setminus\Ga$
with integer coefficients which are all positive or negative.
Thus, the subset $\bO^+\subset\bO$ of the elements with positive
coefficients is exactly the image of $\Omega^+$.
We call elements of $\bO$ quasiroots and the images  of $\Pi\setminus\Ga$
simple quasiroots.

\begin{propn}\label{pr01}
a) Let $\beta$ and $\beta^\prime$ be roots from $\Omega$ such that they give
the same element in $\bO$.
Then there exist roots $\aa_1,\dots,\aa_k\in \Omega_\Ga$
such that $\beta+\aa_1+\cdots+\aa_k=\beta^\prime$ and all partial sums
$\beta+\aa_1+\cdots+\aa_i$, $i=1,\dots,k$ are roots.

b) Let $\bar\beta_1,\dots,\bar\beta_k, \bar\beta$ be elements of
$\bO$ such that $\bar\beta=\bar\beta_1+\cdots+\bar\beta_k$.
Then there exist representatives of these elements in $\Omega$
such that $\beta=\beta_1+\cdots+\beta_k$.
\end{propn}

\begin{proof} a) Let $\beta^\prime=\beta+\gamma_1+\cdots+\gamma_m$ where
$\gamma_i\in \Ga\cup-\Ga$. If $(\beta^\prime,\beta)>0$, then
$\beta^\prime-\beta$ is a root, and the proposition follows.
Proceed by induction on $m$. Since $(\beta^\prime,\beta^\prime)>0$,
if $(\beta^\prime,\beta)\leq 0$, 
then there exists a $\gamma_i$, say  $\gamma_m$, 
such that $(\beta^\prime,\gamma_{m})>0$, so $\beta^\prime-\gamma_{m}$ 
is a root
and $\beta^\prime-\gamma_{m}=\beta+\gamma_1+\cdots+\gamma_{m-1}$.
By induction, the proposition holds for the pair
$\beta^\prime-\gamma_{m}$ and $\beta$, i.e. there exist a representation
$\beta+\aa_1+\cdots+\aa_{k-1}=\beta^\prime-\gamma_m$
satisfying the proposition. Now, putting $\aa_k=\gamma_m$, we obtain
the required representation of $\beta^\prime$.  

b) Let $\pb_1,\dots,\pb_k,\pb$ be some representatives of 
elements $\bbe_1,\dots,\bbe_k,\bbe$ in $\bar\Omega$.
Then we have the equation
$\pb_1+\cdots+\pb_k+\gamma_1+\cdots+\gamma_m=\pb$
for some $\gamma_i\in \Ga\cup-\Ga$, $i=1,\dots,m$.
If $(\gamma_i,\pb)>0$, then $\pb-\gamma_i$ is a root, and
we can take $\pb-\gamma_i$ instead $\pb$. Using iteration, we can regard
that all $(\gamma_i,\pb)\leq 0$. Then, there exists a $\pb_i$, say $\pb_k$,
such that $(\pb_k,\pb)>0$, so $\pb-\pb_k$ is a root.
Applying induction on $k$, one can suppose that there are 
representatives $\beta_1,\dots,\beta_{k-1}$ of $\bbe_1,\dots,\bbe_{k-1},$
such that we have an equation
$\beta_1+\cdots+\beta_{k-1}=\pb-\pb_k+\gamma_1+\cdots+\gamma_n$
for some $\gamma_i\in \Ga\cup-\Ga$, $i=1,\dots,n$,
and $\tilde\beta=\beta_1+\cdots+\beta_{k-1}$ is a root.
If $(\tilde\beta,\gamma_i)>0$, then for some $\beta_j$, say $\beta_1$,
$(\beta_1,\gamma_i)>0$, $\beta_1-\gamma_i$ is a root, and one can
replace $\beta_1$ by $\beta_1-\gamma_i$. Repeating this argument we can
assume that all $(\tilde\beta,\gamma_i)\leq 0$.  
Then $(\tilde\beta,(\pb-\pb_k))\geq 0$, so either
$(\tilde\beta,\pb)>0$, and we set $\beta_k=\pb_k-\sum \gamma_i$,
$\beta=\pb$, or $(\tilde\beta,-\pb_k)>0$, and we set
$\beta_k=\pb_k$ and $\beta=\pb+\sum \gamma_i$.
In any case, we obtain the required representatives of
$\bbe_1,\dots,\bbe_k,\bbe$.
\end{proof}

\begin{rem}\label{rem1}
It is obvious that $\mm$ considered as
a $\g_\Ga$ representation space decomposes into the direct sum of
subrepresentations $\mm_\bbe$, $\bbe\in \bO$,
where $\mm_\bbe$ is generated by all the elements $E_\beta$, $\beta\in \Omega$,
such that the projection of $\beta$ is equal to $\bbe$.
Part  a) of
Proposition \ref{pr01} shows that all $\mm_\bbe$ are irreducible.
Part b) together with part a) shows that for $\bbe_1,\bbe_2\in \bO$ such
that $\bbe_1+\bbe_2\in \bO$ one has
$[\mm_{\bbe_1},\mm_{\bbe_2}]=\mm_{\bbe_1+\bbe_2}$.
Using the Killing form, it is easy to see that representations
$\mm_\bbe$ and $\mm_{-\bbe}$ are dual.

{\it Question}. Is it true that for $\bbe_1,\bbe_2\in \bO$ such
that $\bbe_1+\bbe_2\in \bO$ the representation $\mm_{\bbe_1+\bbe_2}$
is contained in $\mm_{\bbe_1}\wedge\mm_{\bbe_2}\subset \wedge^2 \mm$
 with multiplicity one?
\end{rem}
 
Since $\g_\Ga$ contains the Cartan subalgebra $\hh$, each $\g_\Ga$ invariant
tensor over $\mm$ has to be of weight zero. It follows that there are
no invariant vectors in $\mm$. Hence, there are no invariant vector
fields on $M$.

Consider the invariant bivector fields on $M$.
{}From the above, such fields correspond to the $\g_\Ga$ invariant
bivectors from $\wedge^2\mm$. Note, that any $\hh$ invariant bivector
from $\wedge^2\mm$ has to be of the form 
$\sum c(\aa)\Ea\wedge\Ema$.

\begin{propn}\label{pr0}
A bivector  $v\in\wedge^2\mm$ is $\g_\Ga$ invariant if and only if
it has the form
$v=\sum c(\aa)\Ea\wedge\Ema$
where the sum runs over $\aa\in \rc$,
and for two roots $\aa, \beta$ which give the same element
in $\hh^*/\hh^*_\Ga$ 
one has $c(\aa)=c(\beta)$.
\end{propn}

\begin{proof} 
In view of Proposition \ref{pr01} a), we may assume 
that $\aa=\beta+\gamma$ where $\gamma\in \Omega_\Ga$.
Then the coefficient before 
$E_{\beta+\gamma}\wedge E_{-\beta}$ in $\[E_\gamma,v\]$
appears from the terms
$E_{\beta}\wedge E_{-\beta}$ and $E_{\beta+\gamma}\wedge E_{-\beta-\gamma}$
in $v$, and is equal to
$N_{\gamma,\beta}c(\beta)+N_{\gamma,-\beta-\gamma}c(\beta+\gamma)$.
But from (\ref{he}) follows that $N_{\gamma,\beta}=-N_{\gamma,-\beta-\gamma}$, 
so if $v$ is invariant under action of $E_\gamma$,
i.e. $\[E_\gamma,v\]=0$, then $c(\beta)=c(\beta+\gamma)$.
\end{proof}

This proposition shows, that coefficients of an invariant element
$v=\sum c(\aa)\Ea\wedge\Ema$ depend only of the image of $\aa$ in $\bO^+$,
denoted $\baa$, so $v$ can be written in the form
$v=\sum c(\baa)\Ea\wedge\Ema$.Let $v\in\wedge^2\mm$ be of the form $v=\sum c(\baa)\Ea\wedge\Ema$
where the sum runs over $\aa\in \rc$. Denote by $\theta$ the Cartan
automorphism of $\g$. Then, $v$ is $\theta$ anti-invariant, i.e. 
$\theta v=-v$. Hence, any $\g_\Ga$ invariant bivector is $\theta$ 
anti-invariant.
If  $v,w\in\wedge^2\mm$ are $\g_\Ga$ invariant, then 
$\[v,w\]$ is $\theta$ invariant and is of the form 
$\[v,w\]=\sum e(\baa,\bbe)E_{\aa+\beta}\wedge E_{-\aa}\wedge E_{-\beta}$
where roots $\aa,\beta$ are both negative or both positive
and $e(\baa,\bbe)=-e(-\baa,-\bbe)$.
Hence, to calculate  $\[v,w\]$ for such $v$ and $w$ 
it is sufficient to calculate
coefficients $e(\baa,\bbe)$ for positive $\baa$ and $\bbe$.

\begin{propn}\label{pr1}
Let $v=\sum c(\aa)\Ea\wedge\Ema$, $w=\sum d(\aa)\Ea\wedge\Ema$
be elements from $\g\wedge\g$.\\
Then for any positive roots $\aa,\beta,(\aa+\beta)$
the coefficient by the term 
$E_{\aa+\beta}\wedge E_{-\aa}\wedge E_{-\beta}$ in
$\[v,w\]$
is equal to
\be{}\label{coef}
N_{\aa,\beta}&(c(\aa)(d(\beta)-d(\aa+\beta))+ \\
&c(\beta)(d(\aa)-d(\aa+\beta))-  \notag \\
&c(\aa+\beta)(d(\aa)+d(\beta))).   \notag
\ee{}

\end{propn}

\begin{proof} 
Direct computation, see \cite{KRR}.
\end{proof}

Let $r\in\g\wedge\g$ be the Drinfeld-Jimbo $r$-matrix:
\be{}\label{rm}
r=\sum_{\aa\in\Omega^+}E_\aa\wedge E_{-\aa}.
\ee{}
Then $\[r,r\]=\ff$ is an invariant element in $\wedge^3\g$.
{}From Proposition \ref{pr0} follows that  $r$ reduced  modulo
 $\g \wedge \g_\Ga$ is $\g_\Ga$-invariant. 
Hence, $r$ and $\ff$ define invariant
bivector and three-vector fields on $M$, $\pi^l(r)$ and $\pi^l(\ff)$,
which we denote by $r_M$ and $\ff_M$. Recall, that we identify
invariant tensor fields on $M$ with invariant tensors in $\mm$.

{}From Propositions \ref{pr1} and \ref{pr01} b) 
follows that the condition that the Schouten bracket of
bivector $v=\sum c(\baa)E_{\aa}\wedge E_{-\aa}$ 
with itself give $K^2\ff_M$ for a number $K$ is
\be{}\label{ff} 
c(\baa+\bbe)(c(\baa)+c(\bbe))=c(\baa)c(\bbe)+K^2
\ee{}
for all the pairs of positive quasiroots $\baa, \bbe$ such that 
$\baa+\bbe$ is a quasiroot.

Given $c(\baa)$ and $c(\bbe)$ and assuming that  
$c(\baa)+c(\bbe)\neq 0$ we find that
\be{}\label{ab}
c(\baa+\bbe)=\frac{c(\baa)c(\bbe)+K^2}{c(\baa)+c(\bbe)}.
\ee{}

 Assume, $\baa,\bbe,\bga$ are positive quasiroots such that
$\baa+\bbe, \bbe+\bga, \baa+\bbe+\bga$ are also quasiroots.
Then the number $c(\baa+\bbe+\bga)$ can be calculated formally 
(ignoring possible division by zero) in two ways,
using (\ref{ab}) for  the pair $c(\baa), c(\bbe+\bga)$ on the right hand side
and  also for the pair $c(\baa+\bbe), c(\bga)$.
But it is easy to check that these two ways
give the same value of $c(\baa+\bbe+\bga)$. In this sense the system
of equations corresponding to  (\ref{ab}) for all pairs is consistent.

Let us consider this system more carefully. For any positive quasiroot 
represented as sum of simple quasiroots, 
$\baa=\sum a_i\baa_i$, define the height
 ${\rm ht}(\baa)=\sum a_i$.  In general the coefficient $c(\baa)$ 
 for quasiroots of height $l$ can be formally  defined by iterating 
(\ref{ab}). Let  $\baa=\bar{\alpha}_1+\cdots+ \bar{\alpha}_l$ (with 
possible repetitions) and set $c_i:=c(\bar{\alpha}_i)$, then
\be{}\label{abk}
c(\baa)=\frac{c_{1}c_{2}\cdots c_{l}+
K^2\sum c_{i_1}c_{i_2}\cdots c_{i_{l-2}}+\cdots}
{\sum c_{i_1}c_{i_2}\cdots c_{i_{l-1}}+K^2\sum c_{i_1}c_{i_2}\cdots c_{i_{l-3}}
+\cdots}.
\ee{}
This  expression  can be written as
\be{*}
K\frac{\left(\Pi_{i=1}^l(c_i +K)+\Pi_{i=1}^l(c_i -K)\right)}
{\left(\Pi_{i=1}^l(c_i +K)-\Pi_{i=1}^l(c_i -K)\right)}.
\ee{*}

{\it Assumption}
Let  $(\baa_1,\cdots,\baa_k)$ be the $k$-tuple of simple quasiroots.
In the following we will assume that the point 
$(c_1,\ldots, c_k):=(c(\baa_1),\ldots,c(\baa_k))\in\C^k$ 
does not lie on any of the subvarieties defined by the expressions
in  the denominator of (\ref{abk}) and $c_i\neq 0$ for all $i=1,\dots,k$.

\begin{propn}\label{pr2}
a) Given a $k$-tuple of  positive numbers $(c_1,\ldots,c_k)$ 
satisfying the Assumption, equation
(\ref{abk}) uniquely defines numbers $c(\baa)$ for all 
positive quasiroots $\baa=\sum \baa_i$ such that the bivector
$v=\sum c(\baa)E_{\alpha}\wedge E_{-\alpha}$ satisfies the condition
$\[v,v\]=K^2\ff_M$.\\
b) When $K=0$, the solution  described in part a)
 defines a Poisson bracket on $M$. If all $c(\baa_i)$ are nonzero
(as in Assumption),
there exists a linear form $\ll\in \hh^*_{\Pi\setminus\Ga}$ such that
\be{}\label{ll}
c(\baa)=\frac{1}{\ll(\baa)}
\ee{}
for all quasiroots $\baa$.
\end{propn}

\begin{proof}
a) Since by assumption the denominator is never zero, equation (\ref{abk})
defines consistently $c(\bbe)$ satisfying (\ref{ab}) for all positive $\bbe$.

b) When $K=0$  (\ref{ff}) becomes
$c(\baa+\bbe)(c(\baa)+c(\bbe))=c(\baa)c(\bbe)$.
Assuming $c(\baa_i)\neq 0$ for the simple quasiroots,
 equation (\ref{ff}) implies $c(\baa)\neq 0$ for all quasiroots, so
setting $\ll(\baa)=1/c(\baa)$, we find
 that equation (\ref{ff}) is equivalent to the equation
$\ll(\baa+\bbe)=\ll(\baa)+\ll(\bbe)$.
Thus  $\ll$ is a linear functional, i.e., an element of 
$\hh_{\Pi\setminus \Ga}^*$, which by the construct 
must be nonzero on all quasiroots.
\end{proof}
 
\begin{rem}\label{rem2}
a) This proposition shows that invariant brackets $v$ on $M$ 
such that $\[v,v\]=K^2\ff_M$ form a $k$-dimensional manifold, $\X_K$,
which equals $\C^k$ minus the subvarieties defined in the {\it Assumption},  
where $k$ is the number of elements from $\Pi\setminus\Ga$.
Further, it is known that $k=\dim H^2(M)$, \cite{Bo}. 
If $K$ is regarded as indeterminate, then
$v$ forms a $k+1$ dimensional manifold, $\X\subset \C^k\times\C$,
(component $\C$ corresponds to $K$).
Submanifold $\X_0$ corresponds to $K=0$, i.e. consists of
Poisson brackets. It is easy to see that all the Poisson brackets
of the type $c(\baa)=1/\ll(\baa)\neq 0$ are nondegenerate. Since $\X$ is 
connected, 
it follows that almost all brackets $v$ (except an algebraic subset 
in $\X$ of lesser dimension) are nondegenerate as well.

b) If $v$ defines a Poisson bracket on $M$, then $M$ is a symplectic manifold
and may be realized as an orbit in $\g^*$ passing through
the element $\ll$ from (\ref{ll}) trivially extended to $\hh^*$,
with the KKS bracket.
\end{rem}

Now we fix a Poisson bracket
$v=\sum (1/\ll(\baa))E_{\aa}\wedge E_{-\aa}$ 
where $\ll$ is a fixed linear form
and
describe the invariant brackets 
$f=\sum c(\baa)E_{\aa}\wedge E_{-\aa}$ which satisfy 
the conditions
\be{}\label{pbf}
&\[f,f\]=K^2\ff_M  \qqquad {\rm for} \quad K\neq 0, \\ 
&\[f,v\]=0. \notag
\ee{}

An ordered pairs of quasiroots $\baa,\bbe$ 
such that $\baa+\bbe$ is a quasiroot as well will be called an
{\it admissible pair}.
Substituting in (\ref{coef}) instead $d(\aa)$ the coefficients of $w$,
we obtain that the condition $\[f,v\]=0$ is equivalent 
to the system of equations for the coefficients of $f$
\be{} \label{comp}
c(\baa)\ll(\baa)^2+c(\bbe)\ll(\bbe)^2=c(\baa+\bbe)\ll(\baa+\bbe)^2
\ee{}
for all admissible pairs $\baa$, $\bbe$.

On the other hand, the condition $\[f,f\]=K^2\ff_M$ 
is equivalent to the system of equations
(\ref{ff}) for all admissible pairs of quasiroots.

Substituting $c(\aa+\beta)$ from (\ref{comp}) in (\ref{ff})
we obtain
\be{*}
(c(\baa)\ll(\baa)^2+c(\bbe)\ll(\bbe)^2)(c(\baa)+c(\bbe))=
c(\baa)c(\bbe)(\ll(\baa)+\ll(\bbe))^2+K^2\ll(\baa+\bbe)^2.
\ee{*}
Cancelling terms and extracting the square root, we obtain
the equation
\be{} \label{eq1}
c(\bbe)\ll(\bbe)=c(\baa)\ll(\baa)\pm K\ll(\baa+\bbe).
\ee{}

Substituting $c(\bbe)\ll(\bbe)$ from (\ref{eq1}) in (\ref{comp}), we obtain
\be{} \label{eq2}
c(\baa+\bbe)\ll(\baa+\bbe)=c(\baa)\ll(\baa)\pm K\ll(\bbe)
\ee{}

So, the conditions (\ref{pbf}) 
on $f$ are equivalent to the system of equations (\ref{eq1}, \ref{eq2})
with the same sign before $K$
for all admissible pairs $\baa$ and $\bbe$ from $\bO^+$.

We say that an ordered triple of positive quasiroots (not necessarily 
different) 
$\baa,\bbe,\bga\in \bO^+$ is an {\it admissible triple}, if $\baa+\bbe$, $\bbe+\bga$,
and $\baa+\bbe+\bga$ are quasiroots, too.

\begin{lemma}\label{lem1}
Let $\baa,\bbe,\bga$ be an admissible
triple of quasiroots. If 
$c(\baa+\bbe)\ll(\baa+\bbe)=c(\baa)\ll(\baa) + K\ll(\bbe)$,
then $c(\bbe+\bga)\ll(\bbe+\bga)=
c(\bbe)\ll(\bbe)+ K\ll(\bga)$, that is, the signs before $K$ in  (\ref{eq2}) 
for the admissible pairs $\aa,\beta$ and $\beta,\gamma$ is the same.
\end{lemma}

\begin{proof}
The admissible pair $\baa, \bbe+\bga$ gives the equation
$$c(\baa)\ll(\baa)=c(\bbe+\bga)\ll(\bbe+\bga)\pm K\ll(\baa+\bbe+\bga).$$
The first equation given in the lemma implies that the  $+$ sign  
appears in  equation (\ref{eq1}). These two equations imply
$$c(\bbe)\ll(\bbe)-K\ll(\baa+\bbe)=c(\bbe+\bga)\ll(\bbe+\bga)
\pm K\ll(\baa+\bbe+\bga).$$
Substituting for $c(\bbe+\bga)\ll(\bbe+\bga)$ using (\ref{eq2}), we get
$$c(\bbe)\ll(\bbe)-K\ll(\baa+\bbe)=
c(\bbe)\ll(\bbe)\pm K\ll(\bga)\pm K\ll(\baa+\bbe+\bga),$$
where the last two $\pm$ are independent. However if
$$c(\bbe+\bga)\ll(\bbe+\bga)=c(\bbe)\ll(\bbe)-K\ll(\bga),$$
we have a contradiction, since either  the sign in front of 
$K\ll(\baa+\bbe+\bga)$
is positive and $-K\ll(\baa+\bbe)=+K\ll(\baa+\bbe)$, so $0=\ll(\baa+\bbe)$,
or the sign in front of $K\ll(\baa+\bbe+\bga)$ is negative,
implying  $0=\ll(\bga).$
We conclude that the sign in front of $\ll(\bga)$ must be positive,
\end{proof}
In the situation of the lemma we can express
 $c(\bga)$ in term of $c(\baa)$ as 
\be{}\label{eqg}
c(\bga)\ll(\bga)=c(\baa)\ll(\baa)+K(\ll(\baa+\bbe)+\ll(\bbe+\bga)).
\ee{}  

Now, we consider the pair $(M,\ll)$ as an orbit in $\g^*$ passing through
$\ll$ with the KKS Poisson bracket 
\be{*}
v=v_\lambda=\sum_{\aa\in\Omega^+} \frac{1}{\ll(\baa)}
E_{\aa}\wedge E_{-\aa}.
\ee{*}

\begin{defn}
We call $M$ a {\it good} orbit, if there exists on $M$ an invariant bracket
$f=\sum c(\baa)E_{\aa}\wedge E_{-\aa}$  satisfying 
the conditions (\ref{pbf}).
\end{defn}

\begin{propn}\label{pr5} 
The good semisimple orbits are as follows:

a) For $\g$ of type $A_n$ all semisimple orbits are good.

b) For  all other $\g$, the orbit $M$ is good if and only if
the set $\Pi\setminus\Ga$ consists of one or two roots which appear
in representation of the maximal root with coefficient 1.
\end{propn}

\begin{proof}
a) In this case the system of quasiroots $\bO$ looks like a system
of roots of type $A_k$ for $k$ being equal to the number of elements
of $\Pi\setminus\Ga$. So, the simple quasiroots can be ordered in
a sequence $\bbe_1,\dots,\bbe_k$ in such a way that all subsequences
consisting of three adjacent elements are admissible.
Pick an arbitrary value for $c(\bbe_1)$ and a sign before $K$ in (\ref{eq1})
for the pair $\bbe_1$ and $\bbe_2$. Then, due to Lemma~\ref{lem1}, consistency 
of system (\ref{eq1}, \ref{eq2}) implies that  the sign before $K$ is the same
for all adjacent pairs $\bbe_i$ and $\bbe_{i+1}$.
Using equations (\ref{eq1}) and (\ref{eq2}) for a
fixed sign before $K$ and induction on
$\sht(\aa)$, we find all the coefficients $c(\baa)$ of $f$ and see that
the system (\ref{eq1}, \ref{eq2}) is consistent.

b) Let $\g$ be of type $B_n$ or $C_n$. 
Then the maximal root has the form $\aa_1+2\aa_2+\cdots+2\aa_n$ 
where $\aa_i\in \Pi$. 
Denote $\beta=\aa_2+\cdots+\aa_n$ which is a root.
If $\Pi\setminus\Ga$ does not
contain $\aa_1$, i.e. $\baa_1=0$, then $\bO^+$ contains the admissible pair 
$\bbe,\bbe$. So, from (\ref{eq1}) follows that
$c(\bbe)\ll(\bbe)=c(\bbe)\ll(\bbe)\pm 2K\ll(\bbe)$, i.e.
$\ll(\bbe)=0$ which is impossible.

Assume that $\Pi\setminus\Ga$ contains $\aa_1$ and some roots $\aa_i$ for
$i>1$. Then both $\baa_1$ and $\bbe$ are not equal to zero,
and $\bO^+$ contains the admissible
triple $\bbe,\baa_1,\bbe$. It follows from (\ref{eqg})
that $c(\bbe)\ll(\bbe)=c(\bbe)\ll(\bbe)+2K\ll(\bbe+\baa_1)$, i.e.
$\ll(\bbe+\baa_1)=0$ which is impossible as well, because $\bbe+\baa_1$ is
a quasiroot.

So, for consistency of system (\ref{eq1}, \ref{eq2}) 
in cases $B_n$ and $C_n$, the set $\Ga$ has to contain
all the roots $\aa_i$, $i>1$. But in the latter case the system is
trivially consistent, because in that case the set of quasiroots
looks like $A_1$. The homogeneous space $G/G_\Ga$ is a symmetric space.

Consider the case $D_n$. 

The maximal root has the form
$\aa_1+\aa_2+\aa_3 +2\aa_4+\cdots+2\aa_n$. 
Denote $\beta=\aa_4+\cdots+\aa_n$ which is a root.

Consider several cases. 

The cases when two of $\baa_i$, $i\leq 3$, are equal to zero
and $\bbe$ is not equal to zero lead to an inconsistency in the system
(\ref{eq1}, \ref{eq2}) in the same way as in the cases $B_n$ and $C_n$
considered above.

Assume, that two of $\baa_i$, $i\leq 3$, say $\baa_1,\baa_3$, 
and $\bbe$ are not equal to zero.
Then the sequence 
\be{}\label{adm4}
\baa_1,\baa_2+\bbe,\baa_3,\baa_1+\bbe 
\ee{}
is a sequence of four nonzero quasiroots.
It is easy to see that 
the subsequences 
$\baa_1,\baa_2+\bbe,\baa_3$ and
$\baa_2+\bbe,\baa_3,\baa_1+\bbe$
form admissible triples in $\bO^+$. 
{}From Lemma \ref{lem1} it follows that the sign before $K$ must
be the same in (\ref{eq1}) for all adjacent pairs. Taking, 
for example, the sign plus and applying (\ref{eqg}) 
to the second triple, we obtain the equation
\be{}\label{t1}
c(\baa_1+\bbe)\ll(\baa_1+\bbe)=c(\baa_2+\bbe)\ll(\baa_2+\bbe)+
K(\ll(\baa_2+\bbe+\baa_3)+\ll(\baa_3+\baa_1+\bbe)).
\ee{}
{}From (\ref{eq1}) applied to the pair $\baa_1,\baa_2+\bbe$ 
we obtain
\be{}\label{t2}
c(\baa_2+\bbe)\ll(\baa_2+\bbe)=c(\baa_1)\ll(\baa_1)+
K\ll(\baa_1+\baa_2+\bbe).
\ee{}
Putting $c(\baa_1)$ from (\ref{t2}) in (\ref{t1}) and taking into
account linearity of $\ll$, we obtain the equality
\be{}\label{t3}
c(\baa_1+\bbe)\ll(\baa_1+\bbe)=c(\baa_1)\ll(\baa_1)+
2K\ll(\baa_1+\baa_2+\baa_3+\bbe)+K\ll(\bbe).
\ee{}
On the other hand, expressing $c(\baa_1+\bbe)$ in terms of $c(\baa_1)$
from (\ref{eq2}), we obtain
\be{}\label{t4}
c(\baa_1+\bbe)\ll(\baa_1+\bbe)=c(\baa_1)\ll(\baa_1)\pm K\ll(\bbe).
\ee{}
Now, comparing (\ref{t3}) and (\ref{t4}), we see that if we take
plus before $K$ in (\ref{t4}), then
$\ll(\baa_1+\baa_2+\baa_3+2\bbe)=0$, if we take minus, then
$\ll(\baa_1+\baa_2+\baa_3+\bbe)=0$. But both of the cases are
impossible, since $\ll$ is not equal to zero on quasiroots.

Next, assume that $\bbe=0$ but $\baa_i\neq 0$ for $i=1,2,3$. 
In this case the sequence (\ref{adm4}) makes into the sequence
$\baa_1,\baa_2,\baa_3,\baa_1$. Using the above arguments,
we obtain that in this case must be 
$\ll(\baa_1+\baa_2+\baa_3)=0$ which is impossible, since
$\baa_1+\baa_2+\baa_3$ is a quasiroot.

In the case when $\Pi\setminus\Ga$ contains only one or two roots
of $\aa_i$, $i=1,2,3$, system (\ref{eq1}, \ref{eq2}) is consistent,
because in these cases the set of quasiroots $\bO$ looks like
the system of roots of type $A_1$ or $A_2$.

So, the proposition is proved for the classical $\g$.
For the exceptional $\g$ the proposition follows from the same
arguments.
\end{proof}

\begin{rem}
a) Note that Proposition \ref{pr5} may be reformulated in the following
way:

An orbit $M$ is good if and only if the corresponding system of
quasiroots $\bO$ is isomorphic to a system of roots of type $A_k$.
We say in this case that $M$ is of type $A_k$.

Orbits of type $A_1$ are exactly the orbits which are symmetric spaces. 
For such orbits $\ff_M=0$,
and we may take $f=0$. Symmetric orbits exist for all classical $\g$
and also for $\g$ of types $E_6$, $E_7$.

Orbits of type $A_2$ exist for $\g$ of type $A_n$,
$D_n$, and $E_6$.

Orbits of type $A_k$, $k>2$, exist only in case $\g=sl(n)$.
Moreover, in this case all semisimple orbits have the type $A_k$
for $k\geq 1$.

b) From the proof of Proposition \ref{pr5} it follows that the
bracket $f$ satisfying (\ref{pbf}) is defined on good orbits
by the value of its coefficient $c(\aa)$ for a fixed simple root
and the choice of  a sign before $K$.
On the other hand, if a fixed $f_0$ satisfies (\ref{pbf}), then
the family $\pm f_0+sv$ for arbitrary numbers $s$ also satisfies these
conditions. So, this family consists of all invariant brackets satisfying
(\ref{pbf}). 
Almost all brackets from this family (except a finite number) 
are nondegenerate, since
$v$ is nondegenerate and, therefore, for large $s_0$ the bracket
$\pm f_0+s_0v$ is nondegenerate as well.

For symmetric orbits $\Pi\setminus \Gamma$ consists of one element,
there is one quasiroot and $f_0$ is a multiple of $v$. 
Note, that in \cite{GP} is given a classification of all orbits
in coadjoint representation (not necessarily semisimple) for which
$\ff_M=0$.

In particular, if we take $K=1$ and find $f_0$ such that
$\[f_0,f_0\]=\ff_M$, then the family $\pm if_0+sw$ gives all the brackets
satisfying $\[f_0,f_0\]=-\ff_M$ and compatible with KKS bracket on $M$.

\end{rem}

Note, that if $f$ is a bracket satisfying
$\[f,f\]=-\ff_M$ and $\{\cdot,\cdot\}_r$ is
the $r$-matrix bracket (\ref{rmb}), then $f\pm\{\cdot,\cdot\}_r$ 
is a Poisson bracket on $M$
compatible with KKS bracket.

\section{Cohomologies defined by invariant brackets}

In the next section we prove the existence of a 
$U_h(\g)$ invariant quantiztion 
of the Poisson brackets described above using the methods 
of \cite{DS1}. This requires us to consider the $3$-cohomology of
the complex 
$(\Lambda^\bullet (\g/\g_\Ga))^{\g_\Ga}=(\Lambda^\bullet\mm)^{\g_\Ga}$
of $\g_\Ga$ invariants
with differential given by the Schouten
bracket with the bivector 
$v\in(\Lambda^2\mm)^{\g_\Ga}$ from Proposition 
\ref{pr2} a), 
$$\delta_v:u\mapsto \[v,u\]\,\qqquad\mbox{for} 
\quad u\in(\Lambda^\bullet\mm)^{\g_\Ga}.$$ 
The condition  $\delta_v^2=0$ follows from the Jacobi identity for
the Schouten bracket together with the fact that $\[v,v\]=K^2\ff_M$. The
latter equation is equivalent to 
$\[v,v\]=\ff\mbox{ modulo }\g_\Gamma\wedge\g\wedge\g$,
hence $\[v,v\]$ is  invariant modulo $\g_\Gamma\wedge\g\wedge\g$.
Denote these cohomologies by $H^k(M,\delta_v)$, whereas the usual
de Rham cohomologies are denoted by $H^k(M)$.
 
Recall, Remark 3.2 a), that the brackets $v$ satisfying $\[v,v\]=K^2\ff$ 
form a connected manifold $\X$ which contains a submanifold $\X_0$
of Poisson brackets.

\begin{propn}\label{pr4.1} For almost all $v\in \X$ 
(except an algebraic subset
of lesser dimension) one has
\be{*}
H^k(M,\delta_v)=H^k(M)
\ee{*}
for all $k$.
In particular, $H^k(M,\delta_v)=0$ for odd $k$.
\end{propn}

\begin{proof}
First, let $v$ be a Poisson bracket, i.e. $v\in\X_0$.
Then the complex of polyvector fields on $M$, $\Theta^\bullet$, with
the differential $\delta_v$ is well defined.
Denote by $\Omega^\bullet$ the de Rham complex on $M$.
Since none of the coefficients $c(\baa)$ of $v$ are
zero, $v$ is a nondegenerate bivector field, and therefore
it defines an $\A$-linear isomorphism 
$\tilde{v}:\Omega^1\to\Theta^1$, $\omega\mapsto v(\omega,\cdot)$,
which can be extended up to the isomorphism $\tilde{v}:\Omega^k\to\Theta^k$
of $k$-forms onto $k$-vector fields for all $k$.
Using Jacobi identity for $v$ and invariance of $v$, one can
show that $\tilde{v}$ gives a $G$ invariant isomorphism of
these complexes, so their cohomologies are the same.

Since $\g$ is simple, the subcomplex of 
$\g$ invariants, $(\Omega^\bullet)^\g$, splits off
as a subcomplex of $\Omega^\bullet$.
In addition, $\g$ acts trivially on cohomologies,
since for any $g\in G$ the map $X\to X$,
$x\mapsto gx$, is homotopic to the identity map. ($G$ is assumed connected.)
It follows that cohomologies of complexes
$(\Omega^\bullet)^\g$ and $\Omega^\bullet$ coincide.

But $\tilde{v}$ gives an isomorphism of complexes
$(\Omega^\bullet)^\g$ and 
$(\Theta^\bullet)^\g=((\Lambda^\bullet\mm)^{\g_\Ga},\delta_v)$.
So, cohomologies of the latter complex coincide with de Rham
cohomologies, which
proves the proposition for $v$ being Poisson brackets.

Now, consider the family of complexes
$((\Lambda^\bullet\mm)^{\g_\Ga},\delta_v)$, $v\in\X$. It is clear that
$\delta_v$ depends algebraicly on $v$. It follows from the uppersemicontinuity  
of $\dim H^k(M,\delta_v)$ and the fact that $H^k(M)=0$ for odd $k$,
\cite{Bo}, that $H^k(M,\delta_v)=0$ for odd $k$ and almost all $v\in\X$.
Using the uppersemicontinuity again and the fact that
the number $\sum_k(-1)^k\dim H^k(M,\delta_v)$ is the same for all $v\in\X$,
we conclude that $\dim H^k(M,\delta_v)=\dim H^k(M)$ for even $k$ and
almost all $v$.
\end{proof}  

\begin{rem}
Call $v\in \X$ admissible, if it satisfies Proposition \ref{pr4.1}.
{}From  the proof of the proposition follows that the subset $\D$ 
such that $\X\setminus\D$ consists of admissible brackets
does not intersect with the subset $\X_0$
consisting of Poisson brackets. 

Let $M$ be a good orbit and $f_0+sw$ the family from Remark 3.3 b)
satisfying (\ref{pbf}) for a fixed $K$. Then for almost all numbers $s$ this
bracket is admissible. Indeed, this family is contained in the two parameter
family $tf_0+sw$. By $t=0$, $s\neq 0$ we obtain admissible brackets.
So, there exist $t_0\neq 0$ and $s_0$ such that the bracket
$t_0f_0+s_0w$ is admissible. It follows that the bracket
$f_0+(s_0/t_0)w$ is admissible, too. So, in the family $f_0+sw$
there is an admissible bracket, and we conclude that almost all brackets
in this family (except a finite number) are admissible.

{\it Question}.  Is it true that the set of admissible brackets contains
all the nondegenerate brackets?

\end{rem}

For the proof of existence of two parameter quantization
for the cases $D_n$ and $E_6$ in the next section
we will use the following result on invariant three-vector fields.

Denote by $\theta$ the Cartan automorphism of $\g$.

\begin{lemma}
For either $D_n$ or $E_6$ and one of the
subsets, $\Ga$, of simple roots such that $G_\Ga$ defines a good orbit, 
any  $\g_\Ga$ and $\theta$
invariant element $v$ in $\Lambda^3\mm$ is a multiple of $\ff_M$, that is,
 $$\left(\Lambda^3(\mm\right)^{\g_\Ga}\cong\langle\ff_M\rangle.$$ 
\end{lemma}     

\begin{proof}
Let $\g$ be a simple Lie algebra of type $D_n$ or $E_6$ and
 $\{\aa_1,\ldots,\aa_n\}$ a system of simple roots.
Changing notation slightly from Section 3,
we assume that for $\g=D_n$, $( \aa_i,\aa_{i+1})=-1$
for $i=1,\ldots,n-2$, $(\aa_{n-2},\aa_n)=-1$ with all other
inner products of distinct simple roots are zero, and for $\g=E_6$, 
the non-zero products are $( \aa_i,\aa_{i+1})=-1$ 
for $i=1,2,3,4$ and
$( \aa_3,\aa_6)=-1$. For $\g=D_n$, $\Ga$ is one of the subsets
of simple roots, $\Ga_1=\{\aa_1,\ldots,\aa_{n-2}\},\,
\Ga_2=\{\aa_2,\ldots,\aa_{n-1}\}$, or 
$\Ga_3=\{\aa_2,\ldots,\aa_{n-2}, \aa_{n}\}$. 
For $\g=E_6$, $\Ga=\{\aa_2,\aa_3,\aa_4,\aa_6\}.$ 

The positive quasiroots consist of three elements,
$\baa$, $\baa'$, and $\baa+\baa'$. 
Since a $\theta$ invariant element has the form
$w+\theta w$ for $w\in\mm_{\baa}\t \mm_{\baa'}\t \mm_{-(\baa+\baa')}$,
It is sufficient to show that the space of invariants 
in $\mm_{\baa}\t \mm_{\baa'}\t \mm_{-(\baa+\baa')}$
has dimension
one. We know from Remark 3.1 that the subspaces
 $\mm_{\baa},\mm_{\baa'}, \mm_{\baa+\baa'}$ are
irreducible representations  of $\g_\Ga$
and that $\mm_{\baa}$ and $\mm_{-\baa}$ are 
dual. Therefore the dimension of the space of invariants in 
$\mm_{\baa}\t \mm_{\baa'}\t \mm_{-(\baa+\baa')}$ is the multiplicity of the 
representation $\mm_{\baa+\baa'}$ in the tensor product 
$\mm_{\baa}\t \mm_{\baa'}$.
For $D_n$ and any of $\Ga_i$ the algebra $\g_\Ga\cong A_{n-2}$.
For $\Ga_1,\,\aa=\aa_{n-1}$ and $\aa'=\aa_n$, 
the representations $\mm_\baa$ and $\mm_{\baa'}$ are both isomorphic 
to the dual vector representation for $A_{n-2}$, that is, the contragredient
representation to the representation for the fundamental weight $\ll_{n-2}$,
 $$\mm_{\baa_{n-1}}\cong \mm_{\baa_n}\cong (V^{\ll_{n-2}})^*\cong V^{\ll_1}.$$ 
To see that this is so, note first of all that  $\mm_{\baa_{n-1}}$  
is a lowest weight representation because it has a cyclic vector 
$E_{\aa_{n-1}}$ and  all negative simple root vectors
 of $ \g_\Ga$ annihilate $E_{\aa_{n-1}}$.
The corresponding weight of $A_{n-2}$ is $-\ll_{n-2}$ because
 $(\aa_{n-1},\aa_j)=-(\ll_{n-2},\aa_j)$ if $1\leq j\leq n-2$.
The irreducible lowest weight representations
with lowest weight $-\ll_{n-2}$ is $(V^{\ll_{n-2}})^*$ and 
$A_{n-2}$, $(V^{\ll_j})^*\cong V^{\ll_{n-1-j}}.$
Since the subspaces $\mm_\baa$ and $ \mm_{\baa'}$ of $\g$ have zero
intersection the wedge product  $\mm_\baa\wedge  \mm_{\baa'}$  
projects isomorphically onto the tensor product 
which contains the representation 
$$\mm_{\baa+\baa'}\cong (V^{2\ll_{n-2}})^*\cong V^{2\ll_1}$$
with multiplicity one.

In the cases $\Ga_2$ or $\Ga_3$ the representation $\mm_{\baa_1}$ is the
contragredient representation to the vector representation $(V^{\ll_1})^*
\cong V^{\ll_{n-2}}.$
The representations  $\mm_{\aa'}=\mm_{\aa_n}$ for the case
$\Ga_2$ and $\mm_{\aa'}=\mm_{\aa_{n-1}}$ for $\Ga_3$ are 
 $(V^{\ll_{n-3}})^*\cong V^{\ll_2}$. From the elementary representation theory
of the Lie algebra $sl(n-1)$ (type $A_{n-2}$) we see that
 tensor product $\mm_{\baa}\t \mm_{\baa'}$ contains $\mm_{\baa+\baa'}$ with
multiplicity one.

In the case of $E_6$, $\g_\Ga\cong D_4\cong so(8)$. The representations
$\mm_{\baa},\, \mm_{\baa'}$ and  $\mm_{-(\baa+\baa')}$  
are the three inequivalent irreducible eight dimensional
representations, the two spinor representations and the vector representation.
It is well known that there is a one dimensional space of invariants in
the tensor product $\mm_{\baa}\t \mm_{\baa'}\t \mm_{\baa+\baa'}$. 
\end{proof}

The proof of the lemma gives a positive answer on the {\it Question}
from Remark 3.1 for the particular case.
 
Note that in case of symmetric orbits the space of invariant three-vector
fields is equal to zero.

\section{$\Uh$ invariant quantizations in one and two parameters}

In this section we prove the existence of two types of
$\Uh$ invariant quantization
of the function algebra $\A$ on $M=G/G_\Ga$. 
The first is a one parameter quantization
$$ \mu_h(a,b)=\sum_{n\geq 0} h^n\mu_n(a,b)=
ab +\sum_{n\geq 1} h^n\mu_n(a,b),\quad \mu_1(a,b)=
\frac12 (f(a,b)-\{a,b\}_r),$$ 
as described in Section 2, where $f$ is one of the invariant brackets
found above,  
$$f\in (\Lambda^2 \mm)^\g,\quad \[f,f\]=-\ff_M.$$ 
The second is a two parameter quantization
$$\mu_{t,h}(a,b)=ab+h\mu_1(a,b)+t\mu^\prime_1(a,b)+
\sum_{k,l\geq 1} h^kt^l\mu_{k,l}(a,b).$$
Recall that in this case there are two compatible Poisson brackets 
corresponding to such a quantization:
the bracket $\mu_1(a,b)-\mu_1(b,a)$ is skew-symmetric of the form (\ref{rinv})
and $\mu^\prime_1(a,b)-\mu'_1(b,a)$ is a $U(\g)$ invariant bracket 
$v(a,b),$ which we will assume to be the KKS bracket defined by
identifying $G/G_\Ga$ with an orbit of the coadjoint representation.

We remind the reader of the method in  \cite{DS1}.
The first step is to construct a $U(\g)$ invariant
quantization in the category ${\cA}(U(\g)[[h]],\De,\Phi_h).$
Then we use the equivalence
given by the pair $(Id,F_h)$  between the monoidal categories 
${\cA}(U(\g)[[h]],\De,\Phi_h)$
and ${\cA}(U(\g)[[h]],\wt{\De},{\bf 1})$ 
to define a $\Uh$ invariant quantization,
either $\mu_h F^{-1}_h$ in the one parameter case or 
$\mu_{t,h}F^{-1}_h$ in the two parameter case (see Section 2).
In the first step we used the fact that
$(\Lambda^3 \mm)^{\g_\Ga}=0$ for symmetric spaces.
 In the examples considered in this paper, $(\Lambda^3 \mm)^{\g_\Ga}$
does not necessarily vanish, and we modify the proof
using a method from \cite{DS2} (see also \cite{NV}). 
In the case of  $A_n$ any semisimple orbit is
a good orbit and a different method is required (see \cite{Do} 
where the existence of two parameter quantization 
for maximal orbits is proven).  
For the cases
$B_n$ and $C_n$, the only good orbits are symmetric spaces and 
the quantization was dealt with in \cite{DS1}. 
In  the remaining cases $\g=D_n$ or $E_6$, we proved in Lemma 4.1 that
$(\Lambda^3 \mm)^{\g_\Ga,\theta}=\langle\ff_\A\rangle$,
and a suitable modification of the proof still applies.

\begin{propn}
For almost all (in sense of Proposition 4.1)
$\g$ invariant bracket
satisfying $\[f,f\]=-\ff_M$, there exists a 
multiplication  $\mu_h$ on $\A$ 
$$\mu_h(a,b)=ab +(h/2) f(a,b) +\sum_{n\geq 2}h^n \mu_n(a,b)$$
which  is $\g$ invariant (equation \ref{finv})) and
$\Phi$ associative (equation (\ref{fass})).
\end{propn}
\proof To begin,  consider the multiplication $\mu^{(1)}(a,b)=ab +(h/2) f(a,b)$.
The corresponding  obstruction cocycle  is given by 
$$obs_2=\frac1{h^2}(\mu^{(1)}(\mu^{(1)}\otimes id)-
\mu^{(1)}(id\otimes\mu^{(1)})\Phi)$$
considered modulo terms of order $h$. No $\frac1h$ terms appear because
$f$ is a biderivation and, therefore, a Hochschild cocycle. 
The fact that the presence of
$\Phi$ does not interfere with the cocyle condition and 
that this equation defines a Hochschild $3$-cocycle 
was proven in \cite{DS1}. It is well known that if we restrict to
the subcomplex of cochains given by differential operators, the
differential Hochschild cohomology of $\A$ in dimension $p$ is 
the space of polyvector fields on $M$. 
Since $\g$ is reductive, the subspace of $\g$ invariants splits off
as a subcomplex and has cohomology given by $(\Lambda^p\mm)^{\g_\Ga}$.
The complete antisymmetrization of a $p$-tensor projects
the  space of invariant differential $p$-cocycles onto the subspace
$(\Lambda^p\mm)^{\g_\Ga}$ representing the cohomology. The equation
$\[f,f\]+\ff_M=0$ implies that obstruction cocycle is a coboundary,
and we can find a $2$-cochain $\mu_2$ so that 
$\mu^{(2)}=\mu^{(1)}+ h^2\mu_2$ satisfies 
$$\mu^{(2)}(\mu^{(2)}\otimes id)-
\mu^{(2)}(id\otimes\mu^{(2)})\Phi=0\mbox{  mod }h^2.$$
Assume we have defined the deformation $\mu^{(n)}$ to order $h^n$ 
such that $\Phi$ associativity holds modulo $h^n$, then we define the
$(n+1)^{\mbox{st}}$ obstruction cocycle by
$$obs_{n+1}=\frac1{h^{n+1}}(\mu^{(n)}(\mu^{(n)}\otimes id)-
\mu^{(n)}(id\otimes\mu^{(n)})\Phi)\mbox{ mod }h.$$

In \cite{DS1} (Proposition 4) we showed that the usual 
proof that the obstruction cochain satisfies the cocycle condition
carries through to the $\Phi$ associative case. 
The coboundary of $obs_{n+1}$ appears as the 
$h^{n+1}$ coefficient of the signed sum of the compositions of 
$\mu^{(n+1)}$ with $obs_{n+1}$. 
The fact that  $\Phi=1$ mod $h^2$ together with  the pentagon
identity implies that the sum vanishes
identically, and thus all coefficients vanish, including the coboundary
in question. Let $obs_{n+1}'\in (\Lambda^3\mm)^{\g_\Ga}$ be 
the projection of $obs_{n+1}$ on the totally skew symmetric part, which represents
the cohomology class of the obstruction cocycle. The coefficient of
$h^{n+2}$ in the same signed sum, when projected on the skew symmetric
part is $\[f,obs_{n+1}'\]$ which is the coboundary of $obs'_{n+1}$  
in the complex $(\Lambda^\bullet\mm)^{\g_\Ga},\delta_f= \[f,.\])$. Thus
$obs'_{n+1}$  is a $\delta_f$ cocycle. We have shown in Section 4
that this complex has zero cohomology. Now we  modify 
$\mu^{(n+1)}$ by adding a term $h^n\mu_n$ with $\mu_n\in
(\Lambda^2\mm)^{\g_\Ga}$ and consider the $(n+1)^{\mbox{st}}$ obstruction 
cocycle for $\mu^{\prime(n+1)}=\mu^{(n+1)}+h^n\mu_n$. Since the term we added at
degree $h^n$ is a Hochschild cocyle we do not introduce a $h^n$ term
in the calculation of $\mu^{(n)}(\mu^{(n)}\otimes id)-
\mu^{(n)}(id\otimes\mu^{(n)})\Phi $ and the totally skew symmetric
projection $h^{n+1}$ term has been  modified
by $\[f,\mu_n\]$. By choosing $\mu_n$ appropriately
we can make the $(n+1)^{\mbox {st}}$ obstruction cocycle represent the
 zero cohomology class, and we are able to continue
the recursive construction of the desired deformation.\qed
  
Now we prove the existence of a two parameter deformation 
for good orbits in the cases $D_n$ and $E_6$.
\begin{propn}
Given a pair of $\g$ invariant brackets, $f,v$, on
a good orbit in $D_n$ or $E_6$ 
satisfying $\[f,f\]=-\ff_M$, $\[f,v\]=\[v,v\]=0$, there exists a 
multiplication  $\mu_{h,t}$ on $\A$ 
$$\mu_{t,h}(a,b)=ab+(h/2)f(a,b)+(t/2)v(a,b)+
\sum_{k,l\geq 1} h^kt^l\mu_{k,l}(a,b)$$
which  is $\g$ invariant (equation \ref{finv})) and
$\Phi$ associative (equation (\ref{fass})).
\end{propn}
\proof The existence of a multiplication which is $\Phi$ associative
up to and including $h^2$ terms is nearly identical to the previous 
proof. Both $f$ and $v$ are anti-invariant under the Cartan involution
$\theta$. 
We shall look for a multiplication  $\mu_{t,h}$ such that
$\mu_{k,l}$ is $\theta$ anti-invariant and skew-symmetric for odd $k+l$
and $\theta$ invariant and symmetric for even $k+1$.

 So suppose we have a multiplication 
defined to order $n$,
$$\mu_{t,h}(a,b)=ab+h\mu_1(a,b)+t\mu^\prime_1(a,b)+
\sum_{k+l\leq n} h^kt^l\mu_{k,l}(a,b),$$
with mentioned above invariance properties
and $\Phi$ associative to order $h^n$.

Using properties c) and d) for $\Phi$ from Proposition \ref{p1.1},
direct computation shows that the obstruction cochain, 
$$obs_{n+1}=\sum_{k=0,\ldots, n+1} h^k t^{n+1-k}\beta_k,$$
has the following invariance properties:
For odd $n$, $obs_{n+1}$ is $\theta$ invariant and 
$obs_{n+1}(a,b,c)=-obs_{n+1}(c,b,a)$, and for even $n$,  
and $obs_{n+1}$ is $\theta$ anti-invariant and 
$obs_{n+1}(a,b,c)=obs_{n+1}(c,b,a)$.

Hence, the projection of $obs_{n+1}$ on $(\Lambda^3\mm)^{\g_\Ga}$ is
equal to zero for even $n$. It follows that all the $\beta_k$ are 
Hochschild coboundaries, and the standard argument implies
that the multiplication can be extended up to order $n+1$ with
the required properties. 

For odd $n$, Lemma 4.1 shows that
the projection on $(\Lambda^3\mm)^{\g_\Ga}$ has the form
$$obs_{n+1}=\left(\sum_{k=0,\ldots, n+1} a_kh^k t^{n+1-k}\right)\ff_M.$$
The KKS bracket is given by the
two-vector
$$v=\sum _{\alpha\in\Omega^+\setminus\Omega_\Ga}\frac1{\lambda(\bar\alpha)}
E_\alpha\wedge E_{-\alpha}.$$
Setting 
$$w=\sum _{\alpha\in\Omega^+\setminus\Omega_\Ga}\lambda(\bar\alpha)
E_\alpha\wedge E_{-\alpha},$$
gives
$$\[ v,w\]=-3\ff_M.$$
Defining 
$$\mu^{\prime(n)}=\mu^{(n)}+ \frac{a_0}3 t^n w,$$ 
the new obstruction cohomology class is
$$obs'_{n+1}=(\sum_{k=1,\ldots, n+1} a_kh^k t^{n+1-k})\ff_M.$$
Finally we define 
$$\mu^{\prime\prime(n)}=\mu^{\prime(n)}+\sum_{k=1,\ldots, n} a_kh^{k-1} t^{n+1-k})f$$ 
and get an obstruction cocycle which is zero in cohomology.
Now the standard argument implies that the deformation can be extended
to give a $\Phi$ associative invariant multiplication with the required
properties
of order $n+1$.

So, we are able to continue the recursive construction of the
desired multiplication. \qed 

\begin{rem}
Using the $\Phi_h$ associative multiplications $\mu_h$ and $\mu_{t,h}$ 
from Propositions 5.1 and 5.2
and the equivalence between the monoidal categories 
${\cA}(U(\g)[[h]],\De,\Phi_h)$
and ${\cA}(U(\g)[[h]],\wt{\De},{\bf 1})$ 
given by the pair $(Id,F_h)$ (see Section 2), 
one can define $\Uh$ invariant multiplications,
either $\mu_h F^{-1}_h$ in the one parameter case or 
$\mu_{t,h}F^{-1}_h$ in the two parameter case.
\end{rem}

\end{document}